\documentclass[12pt]{amsart}
\usepackage{amsfonts}
\usepackage[english]{babel}
\usepackage{graphicx}
\usepackage{amscd,color}
\usepackage{amsmath}
\usepackage{amssymb}
\usepackage{tikz}
\usepackage{hyperref}
\usepackage{mathtools}

\usepackage{setspace}
\makeatletter
\@namedef{subjclassname@2020}{%
	\textup{2020} Mathematics Subject Classification}
\makeatother

\theoremstyle{plain}

\newtheorem{corollary}{\bf Corollary}

\newtheorem{lemma}{\bf Lemma}

\newtheorem{proposition}{\bf Proposition}
\newtheorem{remark}{Remark}

\newtheorem{theorem}{\bf Theorem}

\numberwithin{equation}{section}

\title[Serrin-type problem in divergence form on Riemannian manifolds]{Serrin-type problem in divergence form on Riemannian manifolds}

\author[M. Batista, M. Santos, A. Silva and J. Sindeaux]{M. Batista$^{1, *}$, M. Santos$^2$, A. Silva$^3$ and J. Sindeaux$^{4}$}

\address{$^{1}$ CPMAT-IM, Universidade Fe\-deral de Alagoas, Macei\'o,  AL, 57072-970, Brazil}
\email{mhbs@mat.ufal.br}

\address{$^{2,~3}$ Departamento de Matem\'{a}tica, Universidade Federal da Para\'{\i}ba, 58.051-900 Jo\~{a}o Pessoa, Para\'{\i}ba, Brazil.}
\email{marcio.santos@academico.ufpb.br}
\email{djacksonalves2016@gmail.com}

\address{$^4$Universidade Regional do Cariri\\
 	63150-000 Campos Sales, Ceará, Brazil}
\email{joyce.sindeaux@urca.br}

\keywords{Overdetermined problem, Soap Bubble, Pohozaev-type identity, Rigidity.}  
\subjclass[2020]{Primary: 35N25, 53C24; Secondary: 35B50, 53B30, 58J05.}

\thanks{$^\ast$ Corresponding author.}

\begin{document}

\begin{abstract}
In this paper, we investigate an overdetermined boundary value problem of divergence type on bounded domains in Riemannian manifolds with non-negative Ricci curvature. Using integral identities and the $P$-function method, we derive geometric inequalities and rigidity results. Under natural conditions on the nonlinearity, we prove that equality implies the domain is isometric to a Euclidean ball, thereby extending classical symmetry results to the Riemannian setting.
\end{abstract}

\maketitle

\section{Introduction} 

The classical Serrin overdetermined problem \cite{serrin} states that the boundary value problem
\begin{align*} 
\left\{
\begin{array}{lll}
\Delta u=-f(u) & \text{in } \Omega \\
u= 0& \text{on } \partial\Omega \\
u_{\nu}=c& \text{on } \partial\Omega, 
\end{array}
\right.  
\end{align*}
admits a solution if and only if $\Omega \subset \mathbb{R}^n$ is a ball and $u$ is a radial function. Here, $f$ is a smooth real-valued function, $\nu$ denotes the outward unit normal vector field on $\partial \Omega$, and $c$ is a positive constant. Serrin provided a complete solution to this problem, which was notably motivated by questions in fluid dynamics. In his work, he employed the method of moving planes, a technique introduced by Alexandrov, to establish the spherical symmetry of the domain.

Another powerful approach to the Serrin problem is known as the Weinberger method or the $P$-function method. This technique involves introducing a suitable $P$-function associated with the solution and applying the classical maximum principle together with the Pohozaev identity for domains in Euclidean space; see \cite{wein}. Following Weinberger's work, numerous contributions have extended overdetermined problems to the setting of Riemannian manifolds; see, for example, \cite{ciraolo, delay, FK, FR, FRS, Roncoroni}.

Our main objective is to study the following overdetermined boundary value problem:
\begin{align} \label{main}
\left\{
\begin{array}{lll}
 \operatorname{div}\Big(\frac{D u}{\sqrt{1+|D u|^2}}\Big)=f(u) & \text{in } \Omega \\
u= 0& \text{on } \partial\Omega \\
u_{\nu}=c>0& \text{on } \partial\Omega, 
\end{array}
\right.  
\end{align}
where $\Omega$ is a $C^2$ bounded open connected domain contained in a Riemannian manifold $M$, $\operatorname{div}$ denotes the divergence operator in $M$, and $\nu$ is the outward unit normal vector field along the boundary $\partial\Omega$ in $M$, and $u$ and $f$ are smooth functions. The particular case where $f \equiv 0$ is closely related to the theory of minimal surfaces. This case was initially studied for bounded domains in $\mathbb{R}^2$ by S. Bernstein \cite{bernestein}, R. Radó \cite{rado}, and J. Douglas \cite{douglas}. An important extension was later provided by H. Jenkins and J. Serrin \cite{jekins}.

The existence of solutions to the overdetermined problem \eqref{main} was further investigated by W.-M. Ni and J. Serrin; see \cite{Ni} and \cite{serrinex}. Since then, this topic has attracted increasing attention in the literature; see, for example, \cite{alvarez2019} and \cite{lopez}. We also note that when $f(u)$ is replaced by a linear term $ku$, the problem is related to the study of capillary surfaces; see \cite{perko}.

We point out that, when \( f \) is a constant, the problem \eqref{main} is closely related to the uniqueness of graphs with constant mean curvature in product manifolds. In this direction, Fragalà, Gazzola, and Kawohl~\cite{FGK} studied the problem \eqref{main} for a general class of elliptic operators on \( M = \mathbb{R}^n \) with \( f(u) = n \). By combining Minkowski's inequality with the Hopf maximum principle, they proved that \( \Omega \) must be a Euclidean ball, provided that \( \partial\Omega \) is star-shaped. Subsequently, Farina and Kawohl~\cite{FK} applied a suitable Pohozaev-type identity to show that if \( u \) is a solution to \eqref{main}, then \( \Omega \) is necessarily a ball.

Classical approaches to the Serrin problem are deeply intertwined with the proof of the celebrated Alexandrov Soap Bubble Theorem, which asserts that the sphere is the unique compact, embedded hypersurface with constant mean curvature in Euclidean space. 

We recall that in Ros’s proof of Alexandrov's theorem via integral inequalities (see \cite{ros}), a key step involves the application of \emph{Heintze–Karcher inequality}, which states the following: given a $n$-dimensional Riemannian manifold $(M, g)$ with non-negative Ricci curvature, and a bounded domain $\Omega \subset M$ whose boundary $\partial\Omega$ has positive mean curvature $H$, one has
\begin{equation*}
\frac{n-1}{n} \int_{\partial\Omega} \frac{1}{H}~dS_x \geq \operatorname{Vol}(\Omega).
\end{equation*}
Moreover, equality holds if and only if $\Omega$ is isometric to a Euclidean ball.

In this paper, we consider the following suitable \( P \)-function:
\[
P = \frac{n}{w} + F(u),
\]
where \( w\coloneqq\sqrt{1 + |D u|^2} \), \(F(u)=\int_0^u f(t)\,dt\), and \( u \) is a solution of \eqref{main}. We shall prove that this function is super-harmonic under suitable assumptions; see Proposition~\ref{pfunction1}. 

By applying the Hopf maximum principle to the \( P \)-function, we establish a Heintze–Karcher-type theorem specifically adapted to our overdetermined problem~\eqref{main}, as presented below.

\begin{theorem}[Heintze–Karcher-Ros-type theorem]\label{HK}
Let \( M \) be a Riemannian manifold with non-negative Ricci curvature, i.e., \( \mathrm{Ric} \geq 0 \). Let \( u \) be a solution of \eqref{main}, where \( f \) is a non-decreasing function and \( f(0) \neq 0 \). If the mean curvature \( \widetilde{H} \) of \( \partial\Omega \) is negative, then
\[
f^2(0)\int_{\partial\Omega} \frac{1}{\widetilde{H}}~dS_x \leq -n f(0) \int_{\Omega} f(u)~dv.
\]
Moreover, equality holds if and only if \( u \) is a radial function and \( \Omega \) is isometric to a Euclidean ball.
\end{theorem}

\begin{remark}
In our computations, we use the outward unit normal vector field, which may result in a sign convention for the mean curvature of $\partial\Omega$ that differs from that found in other literature.
\end{remark}

On the other hand, Magnanini and Poggesi~\cite{poggesi} obtained an interesting integral identity related to the classical Serrin problem, known as a Soap-Bubble type theorem. In recent years, this class of results has attracted considerable attention and has been extended to a wide variety of overdetermined problems in Riemannian manifolds; see, for instance, \cite{araujo, FP, FRS, Gavitone, Magnanini2017, Magnanini2020, Ruan}.

The following result is a Soap-Bubble-type theorem related to our overdetermined problem~\eqref{main}.

\begin{theorem}[Soap Bubble-type theorem]\label{Soap}
Let \( M \) be a Riemannian manifold with non-negative Ricci curvature, i.e., \( \mathrm{Ric} \geq 0 \), and let \( \Omega \subset M \) be a bounded domain. Suppose that $f$ is a non-decreasing function and \( f(0) = n \). If \( u \) is a solution of \eqref{main}, then
\[
\int_{\partial\Omega} (\widetilde{H} - H_0) \left(\frac{u_\nu}{w}\right)^2dS_x \geq 0,
\]
provided that there exists a constant \( H_0 \) satisfying
\[
-\frac{\int_\Omega f_{+}(u)~dv}{|\partial\Omega|} \leq H_0 < 0,
\]
where \( \widetilde{H} \) denotes the mean curvature of \( \partial\Omega \) and \(f_{+}\) is the positive part of $f$. In particular, if \( \widetilde{H} \leq H_0 \) on \( \partial\Omega \), then $u$ is a radial function and \( \Omega \) is isometric to a Euclidean ball.
\end{theorem}

Recall that a vector field \( \Upsilon \in \mathfrak{X}(M) \) is called a \emph{closed conformal vector field} if
\[
D_Y \Upsilon = \varphi Y,
\]
for all \( Y \in \mathfrak{X}(M) \), where \( \varphi \) is a smooth function on \( M \).

An important class of examples arises from warped product manifolds. Indeed, a warped product \( M = I \times_h N \) is the product manifold \( I \times N \) equipped with the metric
\[
g = dt^2 + h^2(t) g_N,
\]
where \( g_N \) is a Riemannian metric on \( N \). Moreover, the vector field \( \Upsilon = h(t) \partial_t \) is a closed conformal vector field with conformal factor \( \varphi = h'(t) \).

Motivated by \cite{jia}, we derive a Pohozaev identity tailored to problem \eqref{main}, and by applying the \( P \)-function method, we obtain the following result.

\begin{theorem}[Rigidity theorem]\label{rigidity}
Let \( M \) be a Riemannian manifold with non-negative Ricci curvature, i.e., \( \mathrm{Ric} \geq 0 \). Suppose there exists a closed conformal vector field \(\Upsilon\) on \( M \) such that \(
\operatorname{div} \Upsilon = n \varphi,
\)
where \(\varphi\) is a smooth function on \( M \) and positive on $\Omega$. If \( u \) is a solution of \eqref{main} with \( f \) non-decreasing, \( f(0) \neq 0 \), and 
\[
\int_\Omega\Big(F(u) - uf(u) - u\Big\langle D(\ln\varphi), \frac{Du}{w}\Big\rangle\Big)\varphi~dv\geq 0,
\]
then $u$ is a radial function and \( \Omega \) is isometric to a Euclidean ball.
\end{theorem}

\textbf{Outline of the paper.} In Section \ref{SBT}, we prove important geometric results, including the Heintze - Karcher - Ros inequality (Theorem~\ref{HK}) and a Soap Bubble-type theorem (Theorem~\ref{Soap}), both adapted to our overdetermined problem \eqref{main}. These results extend classical ideas to the setting of Riemannian manifolds and set the stage for our main findings. In Section \ref{PTE}, we prove the key rigidity theorem (Theorem~\ref{rigidity}) using the Pohozaev type identity and the properties of closed conformal vector fields together with the \( P \)-function method. This shows that under certain conditions, the domain \( \Omega \) must be isometric to a Euclidean ball. 

\section{A soap Bubble type theorem}\label{SBT}

In this section, we focus on the overdetermined problem \eqref{main}, namely,
\begin{align*}
\left\{
\begin{array}{lll}
 \operatorname{div}\left(\dfrac{D u}{\sqrt{1 + |D u|^2}}\right) = f(u) & \text{in } \Omega, \\
u = 0 & \text{on } \partial \Omega, \\
u_{\nu} = c>0 & \text{on } \partial \Omega,
\end{array}
\right.
\end{align*}
where \( u \) is a smooth function, and \( \Omega \) is a bounded domain contained in a Riemannian \( n \)-manifold \( M \) with non-negative Ricci curvature. 

We note that this overdetermined problem is closely related to the geometry of the graph 
\[
\Sigma = \{ (p, u(p)) \in \overline{M} = M \times \mathbb{R} \mid p \in \Omega \}
\]
of the function \( u \) in the product space with product metric. In this context, the unit normal vector field to the graph \( \Sigma^n \) is given by
\(
N = \frac{1}{w} \big( -D u + \partial_t \big),
\)
where \( w = \sqrt{1 + |D u|^2} \).
So, the \textit{angle function} is given by $$\Theta:=\langle N,\partial_t\rangle=\frac{1}{w}.$$
We also define the second fundamental form by 
$A(X)=-\overline{\nabla}_XN,$
where $\overline{\nabla}$ denotes the Riemannian connection of 
$\overline{M}$ and $\nabla$ denotes the connection of $\Sigma.$ The mean curvature is defined by 
$H=\frac{tr A}{n}.$

Now, let us introduce our $P$-function 
$$P=n\Theta+ F(u),$$
where $\Theta=\frac{1}{\sqrt{1+|D u|^2}}=\frac{1}{w}$ and $F(t)=\int_0^t f(s)\, ds$.

Under the previous notation, we have the following result:

\begin{proposition}\label{pfunction1}
Let \( M \) be a Riemannian manifold with non-negative Ricci curvature, i.e., \( \mathrm{Ric} \geq 0 \), and let \( \Omega \subset M \) be a bounded domain. Assume that \( f' \geq 0 \) and \( f(0) \neq 0 \). Then, the function \( P \) attains its minimum on the boundary \( \partial \Omega \), or \( P \) is constant throughout \( \Omega \). Moreover, if \( P \) is constant, then $u$ is a radial function, ${\rm Ric}(Du,\cdot)=0$, and \( \Omega \) is isometric to a Euclidean ball.
\end{proposition}

\begin{proof}
In fact, note that given a solution $u$ of \eqref{main}, we have that the graph of the function $u$ is a hypersurface of dimension $n$ contained in $\mathbb{R}\times M$ such that its mean curvature is given by $H=\frac{f(u)}{n}.$

We recall the classical Jacobi formula, see Proposition 6 of \cite{ADR07},

$$\Delta\Theta=-\big(|A|^2+{\rm Ric}(N^*,N^*)\big)\Theta-n\partial_t(H),$$
where ${\rm Ric}$ denotes the Ricci tensor of $M$ and $N^*$ denotes the orthogonal projection of the normal vector field $N$ on $TM.$ Since $H=\frac{f(u)}{n},$ then $\partial_t(H)=\langle \nabla H, \partial_t\rangle=\frac{1}{n}f'(u)|\nabla u|^2$
and, therefore, 

\begin{equation}\label {Jacobi}
    \Delta\Theta=-\big(|A|^2+{\rm Ric}(N^*,N^*)\big)\Theta - f'(u)|\nabla u|^2.
\end{equation}

On the other hand, we have that

$$\Delta F(u)=f'(u)|\nabla  u|^2+f(u)\Delta u.$$

A straightforward computation show that $\Delta u=nH\Theta$ and, therefore,
\begin{equation}\label{delta gu}
    \Delta F(u)=f'(u)|\nabla u|^2+[f(u)]^2\Theta.
\end{equation}

From \eqref{Jacobi} and \eqref{delta gu} we get that

$$\Delta P=-n\Theta\Big(|A|^2-\frac{[f(u)]^2}{n}+ {\rm Ric}(N^*,N^*)\Big) - (n-1)f'(u)|\nabla u|^2.$$

From Cauchy-Schwartz  $|A|^2-\frac{[f(u)]^2}{n}\geq 0$, non-negativeness of the Ricci curvature, and $f$ non-decreasing, we conclude that 
$$\Delta P\leq 0.$$

According to the Hopf maximum principle (see Theorem 3.5 in \cite{GT98}), the minimum is either attained on the boundary $\partial\Omega$ or $P$ must be constant.
Moreover, if $P$ is constant, then the graph of the function $u$ is a totally umbilical hypersurface, $f$ is a constant function, and ${\rm Ric}(Du, Du)=0$. In particular, from the non-negativeness of the Ricci curvature, $Du$ belongs to the kernel of ${\rm Ric}$. Thus, $A=\frac{f(0)}{n}g,$ and, in particular, $\nabla u$ is a conformal vector field, that is, 
$$\nabla^2u=\frac{\Delta u}{n}g.$$

On the other hand, since $P$ and $f$ are constants we have 
$$n\Theta+uf(0)=c\in\mathbb{R}.$$
Thus,
$$\nabla^2u=\frac{1}{n}f(0)\Theta g=\frac{f(0)(c-uf(0))}{n^2}g.$$

A straightforward calculation shows that
$$\nabla^2(c-uf(0))=-\frac{(f(0))^2}{n}(c-uf(0))g.$$
Finally, since $c-uf(0)$ is constant along the boundary, from Obata type theorem, see Theorem B (II) of \cite{reilly} or Theorem 1.4 of \cite{CLW21}, we conclude that $\Sigma$ is isometric to a spherical cap of sectional curvature $\frac{f(0)}{n}$. In particular, $\partial\Sigma$ is isometric to a sphere.

Thus, since $\partial\Omega = \partial\Sigma$ is isometric to a sphere and the Ricci curvature on $\Omega$ is non-negative, we can apply Theorem 3 of \cite{X97} to conclude that $\Omega$ is a Euclidean ball. Furthermore, from the proof of Theorem B (II) in \cite{reilly}, it follows that the function $u$, as function on $\Sigma$, depends only on the distance to the north pole of the spherical cap. Via the projection $\Pi_1(x,t) = x$, this implies that $u$, viewed as a function on $\Omega$, depends only on the distance to the center of the ball $\Omega$. So, we conclude the desired result.
\end{proof}
\newpage

Next, we introduce a new quantity associated with the main PDE and deduce a sign for it. Moreover, we derive a rigidity result in the case of equality.

\begin{proposition}\label{pfunction2}
    Let \( M \) be a Riemannian manifold with non-negative Ricci curvature, i.e., \( \mathrm{Ric} \geq 0 \) and let $\Omega\subset M$ be a bounded domain. Suppose that $f'\geq 0$ and $f(0)\neq 0$. If $u$ is a solution of \eqref{main}, we have that 
    $$u_{\nu}\bigg(f(0)+\frac{n\widetilde{H}u_\nu}{\sqrt{1+|D u|^2}}\bigg)\geq 0,$$
    where $\widetilde{H}$ is the mean curvature of $\partial\Omega.$
    In particular, if $\widetilde{H}=-f(0)\frac{\sqrt{1+|D u|^2}}{nu_\nu},$ then $u$ is a radial function and $\Omega$ is isometric to a Euclidean ball.
\end{proposition}
\begin{proof}
      The first equation of \eqref{main} can be rewritten pointwise on $\partial\Omega$ as
    \begin{align}\label{div2}
       \left(\frac{1}{\sqrt{1+|D u|^2}}-\frac{u_\nu ^2}{(\sqrt{1+|D u|^2})^3}\right)u_{\nu\nu}-\frac{(n-1)u_\nu \widetilde{H}(x)}{\sqrt{1+|D u|^2}}=f(0).
    \end{align}
    In fact, on $\partial\Omega$ we have that
    \begin{align}\label{div1}
       f(0)&= \operatorname{div}\left(\frac{D u}{\sqrt{1+|D u|^2}}\right)\nonumber\\
       &=\frac{1}{\sqrt{1+|D u|^2}}\Delta u+\left\langle D\left(\frac{1}{\sqrt{1+|D u|^2}}\right),D u\right\rangle\nonumber\\
       &=\frac{1}{\sqrt{1+|D u|^2}}\Delta u -\frac{u_\nu^2}{(\sqrt{1+|D u|^2})^3}u_{\nu\nu}.
    \end{align}
    and the Laplacian of $u$ in $\partial\Omega$ is given by
    \begin{align*}
       \Delta u &= \sum_{i=1}^{n}\langle D_{e_i}Du,e_i\rangle = \langle D_{\nu}D u,\nu\rangle +\sum_{i=1}^{n-1}\langle D_{e_i}D u,e_i\rangle\\
        &=u_{\nu\nu}+u_\nu\sum_{i=1}^{n-1}\langle D_{e_i}\nu,e_i\rangle =u_{\nu\nu}-u_\nu\sum_{i=1}^{n-1}\langle A_\nu(e_i),e_i\rangle\\
        &=u_{\nu\nu}-u_\nu(n-1)\widetilde{H}(x).
    \end{align*}
  Plugging the previous identity in \eqref{div1} we get \eqref{div2} as desired.

    Now, observe that the derivative of $P$ on the boundary of $\Omega$ is given by
\begin{align}\label{Pnu}
    P_\nu&= \langle DP, \nu\rangle=\langle D(n/w)+D F(u),\nu\rangle\nonumber\\
    &=-\frac{nu_\nu}{(\sqrt{1+|D u|^2})^3}D^2 u(\nu, \nu)+\langle f(0)D u,\nu\rangle\nonumber\\
    &=-\frac{nu_\nu}{(\sqrt{1+|D u|^2})^3}u_{\nu\nu}+u_\nu f(0).
\end{align}
Since $\Delta P\leq 0$, from the Hopf maximum principle (see Lemma 3.4 in \cite{GT98}) , we can guarantee that $P_\eta\leq 0$ along $\partial\Omega=\partial\Sigma$, where $\eta$ is the outward unit vector fiel along $\partial\Sigma$ into $\Sigma$. By other hand, a direct computation gives us that $\eta = \frac{1}{\sqrt{1+c^2}}\nu + \frac{c}{\sqrt{1+c^2}}\partial_t$ along $\partial\Omega$, and so $P_\nu =\sqrt{1+c^2}P_\eta\leq 0$. Using \eqref{Pnu} and \eqref{div2}) we get the result. The equality case follows from Proposition \ref{pfunction1}.
\end{proof}

Now, we are able to provide a demonstration of Theorem \ref{HK} as follows. 

\begin{proof}[\textbf{Proof of Theorem \ref{HK}:}] 
    
Let us consider the positive function $\beta=\displaystyle{\frac{n}{w}}$ along the boundary $\partial\Omega.$ Thus,

$$\frac{\big(u_\nu \widetilde{H}\beta+f(0)\big)^2}{\widetilde{H}}=\beta u_{\nu}\big(u_\nu \widetilde{H}\beta+f(0)\big)+u_\nu\beta f(0)+\frac{[f(0)]^2}{\widetilde{H}}.$$

Since $P_\nu=-(n-1)u_{\nu}\big(u_\nu \widetilde{H}\beta+f(0)\big)\leq 0$ and $\widetilde{H}<0$ we conclude that  

$$u_\nu\beta f(0)+\frac{[f(0)]^2}{\widetilde{H}}\leq 0.$$

On the other hand, by Stokes theorem

$$\int_{\partial\Omega} f(0) u_\nu\beta \,dS_x=nf(0)\int_\Omega \operatorname{div}_\Omega\bigg(\frac{D u}{w}\bigg)\,dv=nf(0)\int_\Omega f(u)\, dv.$$

From the above inequalities, we conclude the desired result. Moreover, if the equalities hold, from Proposition \ref{pfunction2}, we deduce that $\Omega$ is isometric to a Euclidean ball. 

\end{proof}

Next, we present the proof of Theorem \ref{Soap}. It is worth noting that, in certain cases, the volume or area element in integrals may be omitted for the sake of simplicity.

\begin{proof}[\textbf{Proof of Theorem \ref{Soap}:}]
    Firstly, note that
\begin{align*}
\int_{\partial\Omega}\widetilde{H}\left(\frac{u_\nu}{w}\right)^2=&H_{0}\int_{\partial\Omega}\left(\frac{u_\nu}{w}\right)^2+\int_{\partial\Omega}(\widetilde{H}-H_{0})\left(\frac{u_\nu}{w}\right)^2\\
 =&H_0\int_{\partial\Omega}\left(\frac{u_{\nu}}{w}+\frac{1}{H_0}\right)^2-2H_0\int_{\partial\Omega}\frac{u_{\nu}}{wH_0}-H_0\int_{\partial\Omega}\frac{1}{H_0^2}+\int_{\partial\Omega}(\widetilde{H}-H_{0})\left(\frac{u_\nu}{w}\right)^2.
\end{align*}

  From previous equation and Proposition \ref{pfunction2} we conclude that
\begin{align}
     0\leq\int_{\partial \Omega}\frac{u_{\nu}}{w}\bigg(f(0)+\frac{n\widetilde{H}u_\nu}{\sqrt{1+|D u|^2}}\bigg)&=f(0)\int_{\partial\Omega}\frac{u_{\nu}}{w}+n\int_{\partial\Omega}\widetilde{H}\left(\frac{u_\nu}{w}\right)^2\nonumber\\
     &=f(0)\int_{\partial\Omega}\frac{u_{\nu}}{w}+nH_0\int_{\partial\Omega}\left(\frac{u_{\nu}}{w}+\frac{1}{H_0}\right)^2-2n\int_{\partial\Omega}\frac{u_{\nu}}{w}-\frac{n}{H_0}|{\partial\Omega}|\nonumber\\
     &+n\int_{\partial\Omega}(\widetilde{H}-H_{0})\left(\frac{u_\nu}{w}\right)^2.\label{soap}  
\end{align}
Since $\displaystyle\operatorname{div}\left(\frac{D u}{\sqrt{1+|D u|^2}}\right)=f(u)$, we get 
$$\int_{\partial\Omega}\frac{u_\nu}{w}=\int_\Omega f(u).$$
Now, recall that
\begin{align*}
   \frac{1}{H_0}\geq -\frac{\displaystyle\int_{\Omega}f(u)}{|\partial\Omega|} \text{ and } f(0)=n,
\end{align*}
and inserting this in \eqref{soap} we get
    $$\int_{\partial\Omega}(\widetilde{H}-H_{0})\left(\frac{u_{\nu}}{w}\right)^2\geq -H_0\int_{\partial\Omega}\left(\frac{u_{\nu}}{w}+\frac{1}{H_0}\right)^2\geq 0$$
as desired. If the equality holds, we can invoke Proposition \ref{pfunction2} to ensure that $\Omega$ is isometric to a Euclidean ball.

\end{proof}

\section{A Pohozaev type equation and an application}\label{PTE}

Our first result in this section is a Pohozaev type identity related to the problem \eqref{main}.

\begin{lemma}\label{Lemma pohozaev}
    Let $M$ be a Riemannian manifold endowed with a closed conformal vector field $\Upsilon$. Let $\Omega\subset M$ be a bounded domain and let $u$ be a solution of the problem \eqref{main}, then  

\[
\int_\Omega P(u)\varphi\, dv + (n-1)\int_\Omega\Big(F(u) - uf(u) - u\Big\langle D(\ln\varphi), \frac{Du}{w}\Big\rangle\Big)\varphi\, dv -\frac{n}{\sqrt{1+c^2}}\int_\Omega\varphi\, dv=0.
\]
\end{lemma}
\begin{proof}
First, using integration by parts we have
\[
\int_\Omega f(u)\langle \Upsilon, Du\rangle = \int_{\partial\Omega}F(u)\langle\Upsilon, \nu\rangle - n\int_\Omega F(u)\varphi,
\]
where we used that ${\rm div}(\Upsilon)=n\varphi$.
On the other hand, after a tedious computation we obtain
\begin{eqnarray*}
\int_\Omega f(u)\langle \Upsilon, Du\rangle &=&\int_\Omega {\rm div}\Big(\frac{Du}{w}\Big)\langle\Upsilon, Du\rangle\\
&=& \int_{\partial\Omega}\langle\Upsilon, Du\rangle\Big\langle \frac{Du}{w}, \nu\Big\rangle - \int_\Omega\frac{1}{w}\Big\langle Du, \varphi Du + D_\Upsilon Du\Big\rangle \\
&=& \frac{nc^2}{\sqrt{1+c^2}}\int_\Omega \varphi - \int_\Omega\frac{\varphi |Du|^2}{w} - \int_\Omega\frac{1}{w}\langle \Upsilon, D_{Du}Du\rangle.
\end{eqnarray*}

Moreover,
\begin{eqnarray*}
\int_\Omega\frac{1}{w}\langle \Upsilon, D_{Du}Du\rangle &=&\int_\Omega div\Big(\frac{|Du|^2}{w}\Upsilon\Big) - \frac{|Du|^2}{w}div(\Upsilon) - \Big\langle D_\Upsilon\Big(\frac{Du}{w}\Big), Du\Big\rangle\\
&=&\frac{nc^2}{\sqrt{1+c^2}}\int_\Omega\varphi - n\int_\Omega\frac{\varphi |Du|^2}{w} + \int_\Omega\Big\langle \Upsilon, D\Big(\frac{1}{w}\Big)\Big\rangle\\
&=& \frac{nc^2}{\sqrt{1+c^2}}\int_\Omega\varphi - n\int_\Omega\frac{\varphi |Du|^2}{w} + \int_\Omega div\Big(\frac{1}{w}\Upsilon\Big) - \frac{1}{w}div(\Upsilon) \\
&=& \frac{nc^2}{\sqrt{1+c^2}}\int_\Omega\varphi - n\int_\Omega\frac{\varphi |Du|^2}{w} + \frac{n}{\sqrt{1+c^2}}\int_\Omega\varphi - \int_\Omega\frac{n}{w}\varphi.
\end{eqnarray*}
Putting such equalities together we get
\[
n\int_\Omega F(u)\varphi + (n-1)\int_\Omega \frac{\varphi |Du|^2}{w} +\int_\Omega \frac{n}{w}\varphi - \frac{n}{\sqrt{1+c^2}}\int_\Omega\varphi=0
\]
i.e.,
\[
n\int_\Omega F(u)\varphi - (n-1)\int_\Omega \Big(\varphi uf(u) + u\Big\langle D\varphi, \frac{Du}{w}\Big\rangle\Big) +\int_\Omega \frac{n}{w}\varphi - \frac{n}{\sqrt{1+c^2}}\int_\Omega\varphi=0,
\]
where we used integration by parts again in the second term. Handling this equality, we obtain the desired result.

\end{proof}

Using the previous Lemma, we are able to provide a demonstration to Theorem \ref{rigidity}.

\begin{proof}[\textbf{Proof of Theorem \ref{rigidity}}]
    By the integral identity of the Pohozaev-type in the previous Lemma and denoting $\Big(F(u) - uf(u) - u\Big\langle D(\ln\varphi), \frac{Du}{w}\Big\rangle\Big)\varphi$ by $\Phi$ we get 
  \[
\int_\Omega P(u)\varphi\, dv +(n-1)\int_\Omega\Phi~dv -\frac{n}{\sqrt{1+c^2}}\int_\Omega\varphi\, dv= 0.
\]

   By Proposition \ref{pfunction1} we know that the $P-$function defined above  satisfies either $$P>\frac{n}{\sqrt{1+c^2}}\text{ in }\Omega\,\,\,\,\,\, \text{ or }\,\,\,\,\,\,P\equiv\frac{n}{\sqrt{1+c^2}}\text{ in }\Omega.$$ Suppose by contradiction that $P>\displaystyle\frac{n}{\sqrt{1+c^2}}\text{ em }\Omega$. Hence, since $\varphi>0$, it follows that

\[
\int_\Omega \Phi~dv< 0,
\]
which contradicts our hypothesis.
    
\end{proof}

We recall that a closed homothetic vector field $\Upsilon$ on $M$ is a closed conformal vector field such that $D_X\Upsilon=cX,$ for some fixed $c\in\mathbb R$ and all $X\in \mathfrak{X}(M).$

As a direct consequence of Theorem \ref{rigidity} we have the following results.

\begin{corollary}
    Let \( M \) be a Riemannian manifold with non-negative Ricci curvature, i.e., \( \mathrm{Ric} \geq 0 \). Suppose there exists a closed homothetic vector field \(\Upsilon\) on \( \Omega \).
    If \( u \) is a solution of 
\begin{align*} \label{cmc}
\left\{
\begin{array}{lll}
 \operatorname{div}\Big(\frac{D u}{\sqrt{1+|D u|^2}}\Big)=n & \text{in } \Omega \\
u= 0& \text{on } \partial\Omega \\
u_{\nu}=c>0& \text{on } \partial\Omega, 
\end{array}
\right.  
\end{align*}

then, $u$ is a radial function, and \( \Omega \) is isometric to a Euclidean ball.
\end{corollary}

\begin{proof}
In fact, since $f$ and $\varphi$ are constant functions (see notation established before), then the result follows directly from Theorem \ref{rigidity}.
\end{proof}


\section*{FUNDING}
The first and second authors were partially supported by the Brazilian National Council for Scientific and Technological Development [Grants: 308440/2021-8, 402563/2023-9  to M.B, 306524/2022-8 to M.S], and by Coordination for the Improvement of Higher Education Personnel [Finance code - 001], Brazil.


\bigskip

\begin{flushleft}
 
 {\bf Ethical Approval:}  Not applicable.\\
 {\bf Competing interests:}  Not applicable. \\
 {\bf Authors' contributions:}    All authors contributed to the study conception and design. All authors performed material preparation, data collection, and analysis. The authors read and approved the final manuscript.\\
{\bf Availability of data and material:}  Not applicable.\\
{\bf Ethical Approval:}  All data generated or analyzed during this study are included in this article.\\
{\bf Consent to participate:}  All authors consent to participate in this work.\\
{\bf Conflict of interest:} The authors declare no conflict of interest. \\
{\bf Consent for publication:}  All authors consent for publication. \\
\end{flushleft}

\bibliographystyle{amsplain}

\end{document}